\documentclass[12pt,fullpage]{article}
\usepackage{amsfonts}

\font\bbigsym=cmmi10 scaled\magstep4 
\font\bbbigsym=cmmi10 scaled\magstep5 
\newcommand{\BIGP}{\lower2pt\hbox{\bbigsym\char'031}} 
\newcommand{\BBIGP}{\lower3pt\hbox{\bbbigsym\char'031}}

\newcommand{\openE}{\hbox{I\kern-.2em E}} 
\newcommand{\openN}{\hbox{I\kern-.2em N}} 
\newcommand\openR{\hbox{I\kern-.2em R}} 
\newcommand\smopenR{{\hbox{\sevenrm I \kern-.5em R}}} 
\newcommand\openH{\hbox{I\kern-.2em H}}

\renewcommand\epsilon{\varepsilon} 
\renewcommand \phi{\varphi}

\renewcommand{\hat}{\widehat} 
\renewcommand\tilde{\widetilde} 
\renewcommand{\bar}{\overline}

\newcommand\implies{\Rightarrow} 
\newcommand\Verti{\Vert_{\infty}}
\newcommand\Vertv{\Vert_v}

\newcommand{\cov}{{\rm cov}}
\newcommand{\var}{{\rm var} }
\newcommand{\Var}{{\rm Var} }

\newcommand{\done}{$ \ \Box$ }
\newcommand{\remm}{{\rm rem}}

\newcommand{\F}{{\cal F}}

\renewcommand{\H}{{\cal H}}

\newcommand{\N}{{\cal N}}
\renewcommand{\P}{{\cal P}}
\newcommand{\Q}{{\cal Q}}

\renewcommand{\l}{{\ell}}

\def\for{\quad {\rm for} \quad} 
\def\tand{\quad {\rm and} \quad}

\def\tiff{\quad{\rm iff} \quad}

\begingroup

\newtheorem{thm}{Theorem}[section]
\newtheorem{prop}{Proposition}[section]
\newtheorem{lem}[thm]{Lemma}
\newtheorem{cor}{Corrolary}[section]
\newtheorem{defn}[thm]{Definition}
\newtheorem{ex}{Example}[section]
\newtheorem{rem}[thm]{Remark}
\newtheorem{con}{Condition}[section] 
\endgroup

\newcommand\bthm {\begin{thm} \rm}
\newcommand\ethm { \end{thm}}
\newcommand\bprop {\begin{prop} \rm}
\newcommand\eprop { \end{prop}}
\newcommand\blem {\begin{lem} \rm}
\newcommand\elem { \end{lem}}
\newcommand\bcor {\begin{cor} \rm}
\newcommand\ecor { \end{cor}}
\newcommand\bdf {\begin{defn} \rm}
\newcommand\edf { \end{defn}}
\newcommand\bex {\begin{ex} \rm}
\newcommand\eex { \end{ex}}
\newcommand\brem {\begin{rem} \rm}
\newcommand\erem { \end{rem}}
\newcommand\bcon {\begin{con} \rm}
\newcommand\econ { \end{con}}

\newcommand\bitem { \begin{itemize}}
\newcommand\eitem { \end{itemize}}

\newcommand\beq {\begin{equation}}
\newcommand\eeq {\end{equation}}

\newcommand\bea{\begin{eqnarray}}
\newcommand\eea{\end{eqnarray}}

\newcommand\beaa{\begin{eqnarray*}}
\newcommand\eeaa{\end{eqnarray*}}

\newcommand{\Proof} {\noindent \it Proof . \rm}

\newcommand\bbib{\begin{thebibliography}{}}
\newcommand\ebib{\end{thebibliography}{}}

\newcommand\bref{\section*{References} 
\begin{list}{0}{\setlength{\rightmargin}{\leftmargin}}}

\newcommand\eref{\end{list}}

\newcommand{\sect}[1]{\section{#1}\setcounter{equation}{0}}
\newcommand{\subsect}[1]{\subsection{#1}}

\title{Information bounds and efficient estimation in
 a class of censored transformation models}
\author{Dorota M. Dabrowska \\
Department of Biostatistics\\
University of California\\
Los Angeles, CA 90095-1772\\
dorota@ucla.edu}

\begin{document}
\maketitle
\pagestyle{myheadings}
\markright{\rm \today}
\parskip=.2in 
\baselineskip=14pt
\voffset=.5in

\newcommand{\muz}{{\mu \;\mbox{a.e.}\; z}}

\begin{center}
Abstract
\end{center}

Transformation models provide a popular tool for regression
analysis of censored failure time data. 
The most common approach towards parameter estimation in these models
is based on nonparametric profile likelihood method.  
Several authors proposed 
also ad hoc  M-estimators of the Euclidean component of the model. These
estimators are  usually simpler to implement  and many of them
have good practical performance. 
In this paper we consider the form of the information bound
for estimation of the Euclidean parameter of the model
and propose a modification of the inefficient M-estimators to
one-step maximum likelihood estimates.

\newpage

\section{Introduction}

The proportional hazard model, originating in Lehmann \cite{r15} and Cox
\cite{r8},
provides the most common tool for regression analysis of failure 
time data.  The simplest version of this model assumes that a 
failure time $T$ and a covariate $Z$ satisfy the interrelationship
\beq
h(T) = - \beta^T Z + \epsilon
\eeq
where $h$ is an unknown continuous increasing function mapping the support of $T$
onto the real line, $\beta$ is an unknown regression coefficient and
$\epsilon$ is an error term independent of the covariate $Z$ and
having extreme value distribution with density $f(x) = e^x \exp [-e^x]$. 
Alternatively, if $\mu$ is the marginal distribution of the 
covariate and  $\bar A(t,z)$ the cumulative hazard function
of the conditional distribution of the  failure time $T$ given  
 $Z =z$, then the  model stipulates that
\beq
\bar A(t,z) = \Gamma(t) e^{\beta^T z} \quad \muz,
\eeq
where  $\beta$ is an unknown regression coefficient and
$\Gamma$ is an unknown continuous increasing 
function, $\Gamma(0) = 0$, mapping the support of the failure time
$T$ onto the positive half line.

During the past two decades
several authors proposed generalizations of (1.1)-(1.2)  to semiparametric
transformation models specifying the interrelationship between
the conditional hazard function $\bar A(t,z)$ and the transformation 
$\Gamma(t)$ as
\beq
\bar A(t,z) = A(\Gamma(t), \theta|z) \quad \muz,
\eeq
where $\{A( \cdot, \theta|z): \theta \in \Theta\}$ is a family
of conditional cumulative hazards dependent on a finite 
dimensional parameter $\theta$. Following Bickel \it et al. \rm \cite{r4}, 
this family is referred to as the "core model".  Common 
 choices of (1.3) include scale regression models with
core model derived from distributions with decreasing hazard rates,
such as the  frailty distributions with finite mean.  
In particular, the proportional odds ratio model has gained much popularity
as a competitor to the proportional hazard model  
\cite{r1, r9,r10, r13,r16,r19,r21}.  Core models
with increasing or non-monotone hazards were considered in \cite{r5,r6,r7,r9,r12}.

In this paper we consider estimation of the parameter $\theta$ based 
on an iid sample of right-censored failure times.
For purposes of analysis of the 
odds ratio model,
Murphy \it et al. \rm \cite{r16}
proposed to  use  nonparametric 
profile likelihood method. 
The approach
taken was similar to the classical proportional 
hazard model. The model (1.3) was extended to 
include all monotone functions. With fixed parameter $\theta$,
an approximate likelihood function for 
the pair $(\theta, \Gamma)$
was  maximized with respect to $\Gamma$ to obtain an estimate
$\Gamma_{n\theta}$ of the unknown transformation. The estimate
$\Gamma_{n\theta}$ was shown to form a step function placing
mass at each uncensored observation, and 
the parameter
$\theta$ was estimated by maximizing the resulting profile likelihood.
Under certain regularity
conditions on the censoring distribution, the authors showed that the 
estimates are consistent, asymptotically Gaussian at rate $\sqrt n$.
They also proposed to estimate the standard errors of the regression
coefficients by numerically twice-differentiating the log-profile likelihood.
The approach was generalized to other transformation models
whose core models have decreasing hazards in \cite{r14, r19,r20}.

In this paper we take a  different approach
towards construction of efficient estimators of the parameter $\theta$.
To motivate it, let us recall \cite{r4}  that 
 in a regular parametric model $\P = \{P_{\theta, \eta}: \theta \in \Theta,
 \eta \in \H\}$, 
the asymptotic variance of any regular estimator of the parameter
$\theta$ satisfies the bound
$$
\var [\sqrt n (\hat \theta - \theta)] \geq [I_{11}(\theta, \eta) -
I_{12}(\theta, \eta) I_{22}(\theta,\eta)^{-1} I_{21}(\theta, \eta)]^{-1},
$$
where 
$$
I(\theta, \eta) = \pmatrix{I_{11}(\theta, \eta), I_{12}(\theta, \eta) \cr
I_{21}(\theta, \eta), I_{22}(\theta, \eta) }
$$
is the Fisher information matrix with entries
\beaa
I_{11}(\theta, \eta) & = & E \dot \l_{\theta}(X, \theta, \eta)^{\otimes 2} \quad
I_{22}(\theta, \eta) = E \dot \l_{\eta}(X, \theta, \eta)^{\otimes 2}, \\
I_{12}(\theta,\eta) & = & E \dot \l_{\theta}(X, \theta, \eta) 
\dot \l_{\eta}(X, \theta, \eta)^T = I_{21}(\theta, \eta)^T .
\eeaa
Here $\dot \l_{\theta}$ and $\dot \l_{\eta}$ represent
score functions corresponding to the two parameters.  Alternatively, 
$$
\var [\sqrt n (\hat \theta - \theta)] \geq  \left 
[E \l^*(X, \theta,\eta)
\l^*(X, \theta,\eta)^T \right]^{-1},
$$
where $\l^*$ is the efficient score function for estimation of 
$\theta$.
The function  $\l^*$ is the (componentwise) projection of 
the vector of scores $\dot \l_{\theta}$ onto the orthocomplement of
the nuisance tangent space $\dot \P_{\eta}$ spanned by all scores of the
nuisance parameter $\eta$. 
To estimate the parameter $\theta$, we may  consider solving
the score equation
\beq
{1 \over n} \sum_{i=1}^n \l^*(X_i, \theta, \hat \eta(\theta)) =o_P(n^{-1/2}),
\eeq
where $\hat \eta(\theta)$ is an estimate of $\eta$ obtained "for each
fixed $\theta$" in the parameter set $\Theta$. 
To be more precise, we assume that  $\P$ is a submodel 
of a larger family of distributions $\Q$, and there exists
a parameter $\tilde \eta :\Q \times \Theta \to \H$, with the 
property
$\tilde \eta(Q, \theta) = \eta$,
whenever 
$Q = P_{\theta, \eta} \in \P$. We require the estimate $\hat \eta(\theta)$
to be  asymptotically unbiased for estimation of 
 $\tilde \eta(Q, \theta)$ in the larger model $\Q$. If the equation 
(1.4) has a consistent solution $\hat \theta$, then $\hat \theta$ is an efficient
estimate of the parameter $\theta$ under additional conditions
on $\hat \eta(\hat \theta)) $ given in 
[4, Ch.7.7].
These amount to the assumption that  
 $ \l^*(\cdot, \theta, \hat \eta(\hat \theta))$ is a consistent
estimator of the efficient score function and the 
bias of this estimator converges
in probability to 0 at a rate faster than $\sqrt n$.

Turning to transformation models, in the case of uncensored data
Bickel \cite{r2}, Bickel and Ritov \cite{r3} and Klaassen \cite{r13}
used invariance of the model with respect to the group of 
increasing transformations to show that the efficient score function
for estimation of $\theta$ is given by a nonlinear rank statistics.
Its form was derived  using Sturm-Liouville theory.
In the case of censored data, several authors verified existence
of $\sqrt n$ estimators of the unknown transformation and used
them to construct ad hoc estimators of the parameter $\theta$.
Whereas these estimators are inefficient, many of them
have good practical performance \cite{r5,r6,r7,r21} and are 
simpler to implement than the profile likelihood method.
They also apply to a wider class of transformation models.

In section 2, we assume the so-called "non-informative censoring
model" and consider core models whose hazard rates are supported on the whole 
positive half-line, and are finite
and positive at $x=0$. Under a certain integrability condition, 
we derive the form of the information bound and efficient
score function for estimation of the parameter $\theta$.
In section 3 we verify the  integrability condition in the special
case of  
the generalized odds ratio  and the
linear hazard regression models.  
In Section 4 we  construct
a class of $Z$ estimators of the parameter $\theta$. To study its properties,
we shall make the  
assumption that the censoring distribution
has support contained in the support of the failure time distribution
 and its upper  point  forms an atom.
Under mild regularity conditions, 
 we show that the proposed $Z$-estimators have an 
asymptotic distribution not depending on the choice of the
estimator of the unknown transformation. As a by-product, we also  show that  
the parameter $\theta$ can be 
efficiently estimated by solving an equation of the form (1.1)
or by means of one-step MLE.

\setcounter{equation}{0}
\section{Information bound} 

\subsection{Martingale identities}

Throughout the paper we assume that the triple  $(X, \delta, Z)$
represents a nonnegative withdrawal time ($X$), a binary withdrawal 
indicator ($\delta = 1$ for failure and $\delta = 0$ for loss-to-follow-up)
and  a vector of covariates ($Z$). The triple $(X, \delta, Z)$ is defined 
on a complete probability space $(\Omega, \F, P)$ and 
 $X = T \wedge T'$ and $\delta = 1(X = T)$, where $T$ and
$T'$ represent failure and censoring times, respectively. 
We assume that $T$ and $T'$ are conditionally independent given $Z$.

We denote by 
\beaa
\bar A(t,Z) & = & \int_0^t P(X \in du, \delta =1| X \geq u, Z), \\
\bar A_c(t,Z) & = & \int_0^t P(X \in du, \delta =0| X \geq u, Z),
\eeaa
the conditional cumulative hazard functions of the failure and 
censoring time
and set
\beaa
N(t) & = & 1(X \leq t, \delta =1), \quad N_{c}(t) = 1(X \leq t, 
\delta =0), \\
\Lambda(t) & = & \int_0^t Y(u) \bar A(du,Z), \quad \Lambda_{c}(t) =
\int_0^t Y(u) \bar A_c(du,Z),
\eeaa
where $
Y(u) = 1(X \geq u)$. 
Then 
$M(t) = N(t) - \Lambda(t)$ and $M_{c}(t) = N(t) - \Lambda_{c}(t)$ 
form mean zero martingales with respect to the self-exciting
filtration $\{
\F_t: t \geq 0\}$, 
$\F_t = \sigma\{1(X \leq t)\delta, 1(X \leq t)X, 1(X \leq t)(1- \delta),
Z\}$.

Let $\tau_0 = \sup \{t: EY(t) > 0 \}$ and let 
$Q$ and $Q_c$ denote the joint subdistribution functions
$$
Q(t,z) = P(X \leq t, \delta = 1, Z \leq z) \tand
 Q_c(t,z) = P(X \leq t, \delta = 0, Z \leq z).
$$
Then the  failure and censoring counting processes satisfy
\beaa
\int_0^{\tau_0} h(u,Z) M(du) \in L_2^0(P) \tiff h \in L_2(Q),\\
\int_0^{\tau_0} h(u,Z) M_c(du) \in L_2^0(P) \tiff h \in L_2(Q_c). 
\eeaa
Assuming that the function $\bar A$ is continuous, the two martingale
processes are  orthogonal. 
Nan, Edmond and Wellner \cite{r17} showed also that 
any function $b \in L_2^0(P)$ can be represented as a sum
$$
b(X, \delta, Z) = \int_0^{\tau_0} R_1[b](u,Z) M(du) + 
\int_0^{\tau_0} R_2[b](u,Z) M_c(du) + E[b|Z],
$$
where $R_1[b] \in L_2(Q)$  and $R_2[b] \in L_2(Q_c)$ are given by 
\beaa
R_1[b](x,z) & = & b(x, 1,z) - E[b(X, \delta,Z)| X \geq x, Z = z], \\
R_2[b](x,z) & = & b(x, 0,z) - E[b(X, \delta,Z)| X \geq x, Z = z]. 
\eeaa
The three terms in this representation are orthogonal.

\subsection{Assumptions and notation}

To derive the form of the information bound and efficient score 
function in the transformation model, we make the following regularity
conditions on the core model.

\bcon
\bitem
\item[{(i)}]    The parameter set $\Theta \subset R^d$ is open, 
and the   parameter $\theta$ is identifiable in the core model: 
$\theta \not= \theta'$ iff  $A(\cdot, \theta,z)
\not \equiv A(\cdot, \theta',z)$  for $\mu$--a.e. $z$.

\item[{(ii)}] For $\mu$ almost all $z$, the function $A(\cdot, \theta,z)$
has a hazard rate $\alpha(\cdot, \theta,z)$ supported on the whole
positive half-line and there exist
constants $0 < m_1 < m_2 < \infty$ such that
$m_1 \leq \alpha(0, \theta,z) \leq m_2$
for $\mu$--a.e. $z$ and all $\theta \in \Theta$.

\item[{(iii)}] The function 
$\l(x, \theta,z) = \log \alpha(x, \theta,z)$ is  
continuously differentiable with respect to both $x$ and $\theta$.

\item[{(iv)}] The family $\{\alpha(x, \theta,z): \theta \in \Theta\}$
forms a regular parametric model, in particular information is
finite and positive definite.

\eitem
\econ

The derivatives of $\l(x,\theta,Z)$ are denoted by
$$
\dot \l(x,\theta,Z) = {\partial \over \partial \theta} \l(x,\theta,Z) \quad\mbox{and}
\quad \l'(x, \theta,Z) = {d \over dx}\l(x, \theta,Z).
$$

Under assumption of transformation model, the true distribution
$P$ of $(X, \delta, Z)$ belongs to the family $\P = \{P_{(\theta, \eta)}:
\theta \in \Theta, \eta \in \H\}$ where $\theta$ is the Euclidean parameter
of interest and $\eta = (\Gamma, \bar A_c, \mu)$ is the nuisance parameter
corresponding to the unknown transformation $\Gamma$, the unknown 
cumulative hazard function $\bar A_c(t;z)$  of the conditional distribution
of the censoring time given the covariate $Z$  and the marginal
distribution $\mu$ of $Z$. We make the following regularity conditions.

\bcon Let $\tau_0 = \sup\{t: P(X \geq t) > 0\}$ and $\tau_F = \sup\{t: P(T \geq t)
> 0\}$.

\bitem
\item[{(i)}] The distribution $\mu$ is nondegenerate. 

\item[{(ii)}] The parameters $\bar A_c$ and
$\mu$ are noninformative on $(\theta, \Gamma)$. 

\item[{(iii)}] $\Gamma$ is an increasing continuous function, $\Gamma(0) = 0$
and $\lim_{t \uparrow \tau_F} \Gamma(t) = \infty$.

\item[{(iv)}] If $\tau_0 < \tau_F$ then $\tau_0$ is an atom of the
marginal survival function of the censoring time.
\eitem
\econ

The condition (iii) refers to the unobserved model corresponding to uncensored
data. If $\tau_0 < \tau_F$, then the condition (iv) serves to ensure that the
the transformation can be estimated
consistently only within the range $[0, \tau_0]$. We do not know whether
there exists a consistent estimate of the parameter $\theta$
if $\tau_0, \tau_0 < \tau_F$, is a continuity point of the censoring
distribution.

Using conditions 2.1-2.2 and the assumption of the conditional
indpendence of failure and censoring times,
the score operators for the four parameters are
\beaa
\dot \l_1[\theta](X, \delta,Z) & = &\int_0^{\tau_0} \dot \l(\Gamma(u),  \theta,Z) M(du), \\
\dot \l_2[g](X, \delta,Z) & = &\int_0^{\tau_0} \left [g(u) + 
\l'(\Gamma(u),  \theta,Z) \left (\int_0^u gd\Gamma \right) \right] M(du), \\
\dot \l_3 [b](X, \delta, Z) &= &\int_0^{\tau_0} b(u,Z) dM_c(u), \\
\dot \l_4 [c](X,\delta,Z)) & = & c(Z),
\eeaa
where 
$c(Z)  \in  L_2^0(\mu)$, $b(u,Z) \in L_2(Q_c)$ and
\beq
g(X) +  \l'(\Gamma(X), \theta,Z) \int_0^X g d\Gamma \in L_2(Q).
\eeq
The tangent spaces for the four parameters are $\dot \P_i = [\dot \l_i], 
i = 1, \ldots, 4$, where $[\alpha]$ denotes the closed linear span
of the set $\alpha$ in $L_2^0(P)$.  The spaces $\dot \P_i, i = 2,3,4$ are mutually orthogonal,
and so  are the spaces $\dot \P_i, i = 1$, and $i= 3,4$.
The nuisance tangent space is $\dot \P_{\eta} = \dot \P_2 + \dot \P_3 + 
\dot \P_4$.

In the case of the proportional hazard model,
\beq
\bar A(t,z)  = \Gamma(t) e^{\theta^T Z},
\eeq
 Sasieni \cite{r18} and
Nan, Edmond and Wellner \cite{r17} showed that the  martingale
oparator
$$
U(f) = \int_0^{\tau_0} \left (f(x,Z) - E [f(X,Z)|X =x, \delta=1] \right) M(dx), 
\quad f \in L_2(Q),
$$ 
satisfies (i)  $U(f) \perp \dot \P_{\eta}$ in $L_2(P)$; 
(ii) 
$\Pi(f| \dot \P_{\eta}^{\perp}) = 
U(R_1(f))$ for $f \in L_2^0(P)$ and
(iii) $\dot \P_{\eta}^{\perp} = \{U(f): f \in L_2(Q)\}$.
Here $\Pi( \cdot| \dot \P_{\eta}^{\perp})$ denotes the projection
onto $\dot \P_{\eta}^{\perp}$, 
the orthocomplement of the nuisance tangent space.  
These properties entail that 
the efficient score function for estimation of the parameter
$\theta$ in the model (2.2) is given by  $U(f)$, $f(x, Z) = Z$. 

To obtain an extension of this result to the transformation model,
we shall use the following  notation. Firstly, we shall find it 
convenient to denote the "true" parameters $(\theta, \Gamma)$ as
$(\theta_0, \Gamma_0)$. Further, let 
\beaa
s[1](u, \Gamma, \theta) & = & EY(u) \alpha(\Gamma(u), \theta, Z)\\
s[f](u, \Gamma,\theta) & = & EY(u)f(u,Z)\alpha(\Gamma(u), \theta,Z), \\
e[f](u, \Gamma,\theta) & = & {s[f] \over s[1]}(u, \Gamma, \theta) \\
\eeaa
and 
$$
\cov[f_1, f_2]  =  e[f_1 f_2^T] - e[f_1] 
e[f_2]^T \;,  \quad
\var[f]  =  \cov[f,f] \;
$$
Note that 
if failure and censoring times are conditionally independent given $Z$,
then under the assumption of the transformation model,
the conditional distribution of $Z$ given $X = t$ and $\delta =1$ has
form 
\beaa
P(Z \in B| X = t, \delta =1) & = & {1 \over s[1](t, \Gamma_0,,\theta_0)}
\int_B G_c(t,z) F(\Gamma_0(t), \theta_0,z)
\alpha(\Gamma_0(t), \theta_0, z) \mu(dz) \\
& = & e[f](t, \Gamma_0, \theta_0),
\eeaa
where $f(t,Z) = 1(Z \in B)$,
$F(x, \theta_0,z) = [\exp - A(x, \theta_0,z)]$ is the survival function
of the core model and $G_c$ is the conditional survival function of the
censoring time.
More generally,
$$
e[f](t, \Gamma_0, \theta_0) = E [f(X, Z)| X = t, \delta = 1]
$$
and similarly,
$\var[f](t, \Gamma_0, \theta_0)$ and 
$\cov_0[f_1, f_2](t, \Gamma_0, \theta_0)$ are conditional variance and 
covariance
operators. From section 2.1, we also have 
$$
f(X,Z) \in L_2(Q) \tiff  E \int Y(u) f^2(u,Z) \alpha(\Gamma_0(u), \theta_0, Z) 
\Gamma_0(du) < \infty.
$$
Using domintated convergence theorem,  
$$
f(X,Z) \in L_2(Q) \tiff s[f^2](u, \Gamma_0,\theta_0) \in L_1(\Gamma_0).
$$
By noting that 
$$
\Gamma_0(t) = \int_0^t {EN(du) \over s[1](u, \Gamma_0,\theta_0)},
$$
the square integrability condition reduces to 
$$
f(X,Z) \in L_2(Q) \tiff  e[f^2](u, \Gamma_0, \theta_0) \in L_1(Q).
$$

Next,
with some abuse of notation, we shall write
$s[\l'](u, \Gamma, \theta), e[\l'](u,\Gamma, \theta)$ and
$\var[\l'](u, \Gamma, \theta)$ whenever $f(X,Z) = \l'(\Gamma(X),\theta,Z)$.
Let
\bea
C(t) & = & 
\int_0^t \left (s[1](u,  \Gamma, \theta) \right)^{-2} EN(du), \nonumber \\
B(t) & = & \int_0^t \var[\l'](u, \Gamma, \theta) EN(du), \nonumber \\
\P(u,t) & = & \exp - \int_u^t s[\l'](u,\Gamma, \theta) C(du), \\
K(t,t') & = & \int_0^{ t \wedge t'} C(du) \P(u,t) 
\P(u,t') \nonumber \;. 
\eea
Finally, let 
$D[f], D[f](t) = D[f](t, \Gamma, \theta),$ denote the solution to the linear Volterra equation
$$
D[f](t) = - \int_0^t s[f](u, \Gamma, \theta) C(du) - \int_0^t 
D[f](u-)
s[\l'](u, \Gamma, \theta) C(du)
$$
The equation has a unique locally bounded solution given by
$$
D[f](t) = - \int_0^t s[f](u, \Gamma, \theta) C(du) \P(u,t).
$$

We shall need the following integrability condition.

\bcon 
Let
$$
\kappa(\tau_0) = \int \int_{0 < u \leq t \leq \tau_0} C(du) \P(u,t)^2 B(du)
$$
and suppose that $\kappa(\tau_0) < \infty$.
\econ

\blem If the  condition 2.3 holds then  
 $K \in L_2(B \otimes B)$ and 
$D[f] \in L_2(B)$, for any function $f(X,Z)$ such that $e[f^2] \in L_1(Q)$
\elem

This lemma can be verified using Cauchy-Schwartz inequality and dominated
convergence theorem. If failure and censoring times
are conditionally independent and the transformation model
is satisfied by $(\theta, \Gamma) = (\theta_0,\Gamma_0)$, then the second part of the lemma holds for any $f(X,Z) \in L_2(Q)$.

\subsection{Information bound}

Unless 
this leads to confusion,  we suppress the dependence of the functions
 $s[f](\cdot, \Gamma,\theta)$,
$e[f](\cdot, \Gamma,\theta)$ and the corresponding variance and
covariance operators on $(\Gamma, \theta)$. 
We write $s_0[f], e_0[f], \var_0[f]$ and $\cov_0[f_1,f_2]$ whenever
the failure and censoring times are conditionally independent given $Z$
and the transformation model is satisfied by $(\Gamma, \theta) =
(\Gamma_0, \theta_0)$. The functions $C_0, B_0, K_0, \P_0$ and $D_0$
are defined analogously.

Further,
for any $\phi \in L_2(B)$ and $f(X,Z)$ such that $e[f^2] \in   L_1(Q)$, set
\beq
\rho[f, \phi](u) = \cov[f, \l'](u) - \phi(u)\var[\l'](u)
\eeq
and let 
\bea
W_f(t,\theta, Z)  & = & 
f(t,Z) - e[f](t) -
  \left (\l'( \Gamma(t), \theta,Z) - 
e[\l'](t) \right) \phi^*(t) \nonumber \\
& - & {1 \over s[1](t)} \int_t^{\tau_0} \P(t,s) \rho[f,\phi^*](s) EN(ds),
\eea 
where $\phi^*$
is the solution to the Fredholm equation
\beq
\phi^*(t) = - D[f](t) - \int_0^{\tau_0} K(t,u) \phi^*(u) B(du) 
+  \int_0^{\tau_0} K(t,u) \cov[f, \l'](u) EN(du).
\eeq
If $\var[\l'] \equiv 0$, then $W_f(t, \theta,Z) =
f(t,Z) - e[f](t)$.  
In addition, setting 
$\psi^* = \phi^*+ D[f]$, the equation (2.6) simplifies to 
\begin{equation}
\psi(t) - \lambda \int_0^{\tau_0} K(t,u) 
\psi(u) B(du) = \eta(t)  ,
\end{equation}
where $\lambda = -1$ and
\beaa
\eta(t)  & = & 
\int_0^{\tau_0} K(t,u) \rho[f, -D[f]] (u)
EN(du) , \\
\rho[f,-D[f]](u) & = & \cov[f, \l'](u) + \var[\l'](u) D[f](u)
\eeaa
If
$\var[\l'] \not\equiv 0$ then
the kernel $K$ is symmetric, 
positive definite and square integrable with respect to $B$.
Therefore it 
 can have only positive eigenvalues. For
$\lambda = -1$, the
equation (2.7)  has a unique solution  given by
\begin{equation}
\psi^*(t)  =  \eta(u) - 
\int_0^{\tau_0} \Delta(t,u, -1) \eta(u) B(du), \\
\end{equation}
where $\Delta(t,u, \lambda)$ 
is the resolvent corresponding to the kernel $K$.
The resolvent equations 
\begin{eqnarray*}
K(t,u) & = & \Delta(t,u, \lambda) - \lambda \int_0^{\tau_0} 
\Delta(t,w, \lambda) B(dw)  K(w,u) \\
 & = & \Delta(t,u, \lambda) - \lambda \int_0^{\tau_0} 
K(t,w) B(dw) \Delta(w,u, \lambda)    ,
\end{eqnarray*}
applied with $\lambda = -1$ imply 
\begin{equation}
\psi^*(t) = \phi^*(t) + D[f](t) = \int_0^{\tau_0} 
\Delta(t,u,-1) \rho[f,-D[f]](u)
EN(du) .
\end{equation}
If $\var [\l'] \not \equiv 0$  but   $\rho[f, -D[f]] \equiv 0$, then 
the solution to this equation is $\psi^* \equiv 0$, or equivalently,
$\phi^* = - D[f]$.

Define
$$
U(f,\theta) = \int_0^{\tau_0} W(f,\theta,Z)(u) M(du)
$$

\bprop
 Suppose that the conditions 2.1-2.3  are satisfied with $(\theta,\Gamma)
= (\theta_0, \Gamma_0)$. 
Then (i)  $U(f,\theta_0) \perp \dot \P_{\eta}$ in $L_2(P)$; 
(ii)
$\Pi(f| \dot \P_{\eta}^{\perp}) = 
U(R_1(f),\theta_0)$ for $f \in L_2^0(P)$ and
(iii) $\dot \P_{\eta}^{\perp} = \{U(f,\theta_0): f \in L_2(Q)\}$.
\eprop

\Proof 
Set
\beaa
g^*(t) & = & e_0[f](t) -
e_0[\l'](t)  \phi^*(t) \\ 
& - & {1 \over s_0[1](t)} \int_t^{\tau_0} \P_0(t,s) \rho_0[f,\phi^*](s) EN(ds).
\eeaa
Then
$$ 
\phi^*(t) = \int_0^t g^* d \Gamma_0
$$
and
$$
U(f,\theta_0) = \int_0^{\tau_0} f(t,Z) M(du) - \dot \l_2[g^*]
= \int_0^{\tau_0} [f(t,Z) 
- g^*(t) - \l'(
\Gamma_0(u), \theta_0,Z) \phi^*(u)] M(du).
$$
Using Cauchy-Schwartz inequality and the condition 2.3, it is easy
to verify that $\dot \l_2[g^*] \in L_2^0(P)$ for any function
$f(X,Z) \in L_2(Q)$.
Therefore part (i) of the proposition will be verified if 
we show that 
$E (U(f,\theta_0) \dot \l_2[g] ) = 0$ for any function $g$ satisfying th
condition (2.1). 

Put
$
\dot \l_2[g] = I_1(g) + I_2(g)$,
where
\beaa
I_1(g) & = & \int_0^{\tau_0} g(u) dM(u), \\
I_2(g) & = & \int_0^{\tau_0} \left [
\l'(\Gamma_0(u), \theta_0, Z) \int_0^u g d \Gamma_0
\right ] M(du). 
\eeaa
Then
\beaa
 E U(f,\theta_0) I_1(g) & = & - \int_0^{\tau_0} g(t) \left [
\int_t^{\tau_0} \P_0(t,s) \rho_0[f,\phi^*](s) EN(ds)
\right ] \Gamma_0(t), \\
 E U(f,\theta_0) I_2(g) 
& = & 
\int_0^{\tau_0} \left [\rho_0[f,\phi^*](t) \int_0^t g d\Gamma_0 \right] dEN(t) \\
& - & \int  
e_0[\l'](t) \left (\int_t^{\tau_0} \P_0(t,s)  
\rho_0[f,\phi^*](s) EN(ds)
\right ) \left (\int_0^t g d\Gamma_0 \right) d\Gamma_0(t).
\eeaa
Using
$$
\P_0(t,s) = 1 - \int_t^s s_0[\l'] dC_0(du)\P_0(u,s) = 1 - \int_t^s  e_0[\l'](u) d\Gamma_0(u) \P_0(u,s)
$$
for $t \leq s$, and applying Fubini theorem, we have
\beaa
& & E U(f,\theta_0)I_2(g)  = 
\int_0^{\tau_0} g(u) \left (\int_u^{\tau_0} 
 \rho_0[f,\phi^*](t) dEN(t) \right ) d \Gamma_0(u) \\
& - & \int_0^{\tau_0} g(u) \left (\int_u^{\tau_0}
e_0[\l'](t) \left (\int_t^{\tau_0} \P_0(t,s)  
\rho_0[f,\phi^*](s) EN(ds)
\right ) d\Gamma_0(t) \right)  d \Gamma_0(u)  \\
& = &
\int_0^{\tau_0} g(u) \left (\int_u^{\tau_0} 
 \rho_0[f,\phi^*](t) dEN(t) \right ) d \Gamma_0(u) \\
& - & \int_0^{\tau_0} g(u) \left (\int_u^{\tau_0} \left (\int_t^s 
e_0[\l'](t) d \Gamma_0(t) \P_0(t,s) \right )  
\rho_0[f,\phi^*](s) EN(ds)
 \right)  d \Gamma_0(u)   \\
& = & \int_0^{\tau_0} g(u) \left (\int_u^{\tau_0} 
 \rho_0[f,\phi^*](t) dEN(t) \right ) d \Gamma_0(u) \\
& + & \int_0^{\tau_0} g(u)  \left (\int_u^{\tau_0} [\P_0(u,s) - 1] 
\rho_0[f,\phi^*](s) EN(ds)
 \right)  d \Gamma_0(u)  = - EU(f,\theta_0) I_1(g).
\eeaa
This completes the proof of part (i). 
Part (ii) and part (iii) follows in the same way as in
Nan, Edmond and Wellner \cite{r17}. \done

Since $U(f, \theta_0)$ is a martingale operator, we have
\beaa
\Sigma_0(f,\theta_0) &= &EU(f,\theta_0)^2 = \\
& = & E \int W_f^2(u,\theta_0,Z) Y(u)\alpha(\Gamma_0(u), \theta_0, Z) \Gamma_0(du). \\
\eeaa          
The matrix $\Sigma_0(f,\theta_0)$ satisfies $\Sigma_0(f,\theta_0) = 
\Sigma_1(f,\theta_0) + \Sigma_2(f,\theta_0)$,
where 
\bea
\Sigma_1(f,\theta_0) & = & \int_0^{\tau_0} \var_0 [f - \l' \phi^*](u) EN(du), \\
\Sigma_2(f,\theta_0) & = & \int_0^{\tau_0} C_0(du) 
\left [\int_u^{\tau_0} \P_0(u,t) \rho_0[f,\phi^*](t) 
EN(du) \right]^{\otimes 2} \;. \nonumber
\eea
The conditional variance function is identically equal to 0, if and only if
\beq
f(t,Z) =
\l'(\Gamma_0(t), \theta_0, Z) \int_0^t h d\Gamma_0 + a(t) \in L_2(Q)
\eeq
 and $a \equiv h$.
For any  function of the form (2.11),  Fubini
theorem yields
\beaa
D_0[f](t) & = & -\int_0^t a(u) d\Gamma_0(u)\P_0(u,t) \\
& - & \int_0^t \left (\int_0^u hd\Gamma_0 \right) e_0[\l'](u) d\Gamma_0(u) \P_0(u,t) \\
& = & - w(t) - \int_0^t hd\Gamma_0, \\
w(t) & = & \int_0^t [a - h](u) d\Gamma_0(u) \P_0(u,t). 
\eeaa
Hence  (2.6) reduces to 
$$
\int_0^t [g^* -h] d \Gamma_0 + 
\int_0^{\tau_0} K_0(t,u) [\left [\int_0^u (g^* -h)d \Gamma_0 \right ] B_0(du) = 
w(t). 
$$
If $a \equiv h$ then the right-hand side of this 
equation is identically equal to 0, and correspondingly, 
 $g^* = h$. Otherwise,
 we obtain
$$
\phi^*(t) = \int_0^t h d\Gamma_0 + w(t) -\int_0^{\tau_0} \Delta(t,u,-1) w(u) B_0(du)
$$
and
$$
\var_0[f - \l' \phi^*](t) = \var_0[\l'](t)[w(t) -  \int_0^{\tau_0} \Delta(t,u,-1) w(u) B_0(du)]
$$
If $\var_0[\l'] \not \equiv 0$, then the right-hand side is 
identically equal to 0 if and only if
$$
w(t) = \int_0^{\tau_0} \Delta(t,u,-1) w(u) B(du).
$$
In this case,  the resolvent equations imply that $w \equiv 0$ so that 
$a \equiv h$.

\sect{Examples}

In this section we verify the square integrability condition
in two models. In particular, we show that the information
bound applied to both censored and uncensored data.

\bex Generalized odds ratio model. 
The survival function of the core model is given by
\beaa
F(x, \theta,Z) & = & [1 + \eta e^{\theta^T Z}x]^{-1 -1/\eta} 
\for \eta > 0, \\
& = & \exp [- e^{\theta^T Z} x] \quad \for  \eta = 0 .
\eeaa
The proportional hazard model corresponds to the choice $\eta = 0$ and
the proportional odds model to $\eta = 1$.
The  hazard function of the core model  given by
$$
\alpha(x, \theta, Z) = e^{\theta^T Z} [1 + \eta e^{\theta^T Z} x]^{-1}.
$$
We have 
\beaa
\dot \l(x, \theta,Z) & = & Z [1 + \eta e^{\theta^T Z} x]^{-1}, \\
\l'(x, \theta, Z) & = & - \eta e^{\theta^T Z} 
[1 + \eta e^{\theta^T Z} x]^{-1}.
\eeaa
If $|Z| \leq d_0$ and $d_1 \leq e^{\theta^T Z} \leq d_2$,
then
\beaa
|\dot \l(x, \theta, Z)| & \leq & d_0 [1 + \eta d_1 x]^{-1} \leq d_0\\
\eta d_1 [1 + \eta d_1 x]^{-1} & \leq & - \l'(x, \theta,Z) \leq \eta d_2
[1 + \eta d_2 x]^{-1}
\eeaa
Note that $-\l'(x, \theta,Z)$ is an increasing function of 
$\eta e^{\theta^T Z}$.

Suppose that $(\theta_0, \Gamma_0)$ is the "true" parameter of the
transformation model. The preceding bounds imply
$$
\P_0(u,t) \leq \exp \int_u^t {\eta d_2 \over 
1 + \eta d_2 \Gamma_0(v) } d\Gamma_0(v) = {1 + \eta d_2 \Gamma_0(t) \over 
1 + \eta d_2 \Gamma_0(u) }.  
$$
Next recall  that if $U$ is a random variable with finite
variance and  distribution function
$H$, then 
\beq
\Var (U) = {1 \over 2} \int \int (w_1 - w_2)^2 H(dw_1) H(dw_2)
\eeq
By noting that  
\beaa 
& & \left [\l'(\Gamma_0(t),\theta_0,z_1) -
\l'(\Gamma_0(t),\theta_0,z_2) \right]^2 = \\
& = & {\eta^2 [e^{\theta^T_0 z_1} - e^{\theta^T_0 z_2}]^2 \over
[1 + \eta e^{\theta^T_0 z_1} \Gamma_0(t)]^2
[1 + \eta e^{\theta^T_0 z_2} \Gamma_0(t)]^2} \\
& \leq & \eta^2 [d_2 - d_1]^2 [1 + \eta d_1 \Gamma_0(t)]^{-4},
\eeaa
and using (3.1), we obtain
$\var [\l'](t, \Gamma_0, \theta_0) = O(1) [1 + \eta d_1 \Gamma_0(t)]^{-4}$.
Since  $[1 + \eta d_2 \Gamma_0(t)] \leq (d_2/d_1) 
[1 + \eta d_1 \Gamma_0(t)]$, we see that  
$$
\kappa(\tau_0) = O(1) \int
\left [ \int_0^t {d \Gamma_0(v) \over 
s[1](v, \Gamma_0, \theta_0)}{ 1 \over [1 + \eta d_2 \Gamma_0(v)]^2} \right] 
{s[1](t, \Gamma_0, \theta_0)  d\Gamma_0(t) \over (1 + \eta d_1 \Gamma_0(t))^{2}} 
$$
In the generalized odds ratio model, $s[1](t, \Gamma_0, \theta_0)$ is a decreasing function of 
$t$. Therefore 
$$
{s[1](t, \Gamma_0, \theta_0) \over s[1](u, \Gamma_0, \theta_0)} \leq 1
$$
and
$$
\kappa(\tau_0) \leq \tilde \kappa(\infty) = O(1) \int [1 - 
(1 + \eta d_2 w))^{-1}] 
(1 + \eta d_1 w)^{-2} dw = O(1)
$$
The bound is valid for any distribution of the censoring times.
We obtain a better bound on the constant $\tilde \kappa(\infty)$ in the case of 
the so-called Koziol-Green censoring
model. The conditional survival function of the censoring time
given the covariate is of the form
$$
G_c(t|z) = F(\Gamma_0(t),\theta_0,z)^{a}, \quad a \geq 0.
$$
The choice of $a = 0$ corresponds to the case of uncensored data.
We have
$$
s[1](t, \Gamma_0, \theta_0) = E e^{\theta_0^T Z} 
[1 + \eta e^{\theta^T_0 Z} \Gamma_0(t)]^{-[ 1 + (1 + a)/\eta]} 
$$
so that 
$$
d_1  [1 + \eta d_2 \Gamma_0(t)]^{-[ 1 + (1 + a)/\eta]}  \leq
s[1](t, \Gamma_0, \theta_0) \leq d_2  [1 + \eta d_1 \Gamma_0(t)]^{-[ 1 + (1 + a)/\eta]}. 
$$
Hence 
\beaa
\tilde \kappa(\infty) & = & O(1) \int_0^{\infty} \left [
\int_0^x (1 + \eta d_2 w)^{(1+a)/\eta -1} dw \right] 
(1 + \eta d_1 x)^{-[3+ (1+a)/\eta]} dx \\
& = & O(1) \int_0^{\infty} (1 + \eta d_1 x)^{-3} dx = O(1). 
\eeaa
\eex

\bex Linear hazard model has failure  rate
$$
\alpha(x, \theta|z) = a_{\theta}(z) + x b_{\theta}(z).
$$
We assume that 
$ m_1 \leq a_{\theta}(z) \leq m_2$, $m_1' \leq b_{\theta}(x) \leq m_2'$
for some finite positive constants $m_q, m'_q, q = 1,2$. In addition,
the functions $a_{\theta}(z)$ and $b_{\theta}(z)$ have bounded derivatives
with respect to $\theta$. We have
$$
\l'(x, \theta,z)  = {b_{\theta}(z) \over 
a_{\theta}(z) + x b_{\theta}(z)} 
$$

Set $d_1 = m_1'/m_2$ and $d_2 = m_2'/m_1$ and suppose that $(\theta_0, \Gamma_0)$
is the true parameter of the transformation model.
Using a similar algebra as in the example 3.1, we can show that 
$$
\P_0(u,t) \leq {1 +  d_1 \Gamma_0(u) \over 
1 +  d_1 \Gamma_0(t) }  \quad 
\var [\l'](t, \Gamma_0, \theta_0) \leq O(1) [1 +  d_1 \Gamma_0(t)]^{-4}.
$$
Hence 
$$
\kappa(\tau_0) = O(1) \int_0^{\tau_0}
\left [ \int_0^t {d \Gamma_0(v) \over 
s[1](v, \Gamma_0, \theta_0)} [1 +  d_1 \Gamma_0(v)]^2 \right] 
{s[1](t, \Gamma_0, \theta_0)  d\Gamma_0(t) \over (1 +  d_1 \Gamma_0(t))^{6}} 
$$
For $v < t$, we have
$$
{s[1](t, \Gamma_0, \theta_0) \over
s[1](v,\Gamma_0, \theta_0)} \leq O(1) {1 + d_2 \Gamma_0(t) \over 1 + d_1 \Gamma_0(v)}
= O(1) {1 + d_1 \Gamma_0(t) \over 1 + d_1 \Gamma_0(v)} 
$$
Hence 
$$
\kappa(\tau_0) \leq  \tilde \kappa(\infty) = O(1) \int_0^{\infty} (1 + d_1 w)^{-3}dw
$$

\eex

Other examples satisfying the integrability condition 2.3 include 
the inverse Gaussian core model, and half-symmetric distributions
such as the half-normal, half-logistic and half-t distribution.

\sect{Estimation}

We turn now to estimation of the parameter $\theta$.
Let us recall that under the assumption of transformation
model, the true distribution of $(X, \delta, Z)$ is in a
class $\P =  \{P_{\theta,\eta}: \theta  \in \Theta,
\eta \in \H\}$, where $\eta$ represents the triple $\eta =( \Gamma, 
G_c, \mu)$.

To construct an estimator of the parameter $\theta$, we assume that
$\Q$ is a class of probability distributions of the variables
$(X, \delta, Z)$ containing $\P$ as a submodel. For each $(Q, \theta)
\in \Q \times \Theta$, we let $\Gamma_{Q, \theta}$ be a monotone
function such that 
$$
\Gamma_{Q, \theta_0} = \Gamma_0
$$
if $Q = P_{\theta_0, \eta_0}$ and $\xi_0 = (\theta_0, \Gamma_0, G_c, \mu)$ 
is the true parameter of the transformation model.

Dropping dependence of this function on the distribution $Q$,
let $\Gamma_{n, \theta}$ be an estimator of  $\Gamma_{\theta}$
such that $\Vert \Gamma_{n\theta} - \Gamma_{\theta} \Verti = o_Q(1)$,
i.e. the estimate is consistent when observations $(X, \delta, Z)$ are 
sampled from a distribution $Q \in \Q$. In addition to this we assume
the following regularity conditions. 

\bcon
Let $B(\theta_0, \epsilon_n)$ denote  an open  ball of radius $\epsilon_n$ 
and centered at $\theta_0$.

\bitem
\item[{(i)}] 
$\epsilon_n \downarrow 0$ and $\sqrt n \epsilon_n \uparrow \infty$.
\item[{(ii)}]  The point $\tau_0 = \sup\{t: EY(t) > 0\}$ is an atom of the
marginal distribution of the censoring times.
\item[{(iii)}] The estimate of the transformation satisfies:
 $\sqrt n \Vert \Gamma_{n0} - \Gamma_0 \Verti = O_P(1)$,
$\limsup_n \{\Vert \Gamma_{n\theta} \Vertv: 
\theta \in B(\theta_0, \epsilon_n)\}  = O_P(1)$ and
$$\sup \{ {\sqrt n \Vert \Gamma_{n\theta} - 
\Gamma_{n0} \Verti/ [\sqrt n |\theta - \theta_0| + 1]}: \theta \not= \theta_0, \theta
\in B(\theta_0, \epsilon_n)\}  = O_P(1)
$$
\eitem
\econ

Examples of estimators satisfying these conditions were given by Cuzick [9],
Bogdanovicius and Nikulin [5] and Yang and Prentice \cite{r21}, among others.

Referring to the notation of section 2, we assume that the function
$f(u, Z)$ is of the form $f(u,Z) = f(\Gamma_{\theta}(u), \theta,Z)$ and
make the following regularity conditions.

\bcon
Let $\psi$ be a constant or a bounded continuous strictly decreasing
function. For $p = 1,2,3$, let $\psi_p$ be continuous bounded or strictly
increasing functions such that $\psi_p(0) < \infty$ and
$$
\int_0^{\infty}  e^{-x} \psi_1^2(x) dx < \infty, \quad  
\int_0^{\infty}  e^{-x}  \psi_2(x) dx < \infty,  \quad
\int_0^{\infty}  e^{-x}  \psi_3(x) dx < \infty \;.
$$
Suppose that the
derivatives of the function $\l(x, \theta,z)$ satisfy
$$
|\l'(x, \theta,z)| \leq \psi(x), \quad  |\l''(x, \theta,z)| \leq \psi(x), \quad
|\dot \l(x, \theta, Z) \leq \psi_1(x)
$$
The function $f(x, \theta,Z)$ is differentiable with
respect to  $x$ and
$$
 |f(x, \theta,Z)| \leq \psi_1(x), \quad \quad |f'(x, \theta,Z)| \leq \psi_2(x).
$$
We also have
\beaa
& & |g_1(x, \theta,z) - g_1(x', \theta,z)| \leq max [\psi_3(x), \psi_3(x')]|x - x'|, \\
& & |g_2(x, \theta,z) - g_2(x, \theta',z)| \leq \psi_3(x) 
\Vert \theta - \theta' \Vert ,
\eeaa
where $g_1 =  \dot \l, f',  \l''$ and $g_2 = \dot \l, f, f', \l', \l''$.   
\econ

The functions $\psi, \psi_p, p = 1,2,3$ may differ in each inequality.
Note that we do not  require differentiability of $f$ with respect to $\theta$,
but only a Lipschitz continuity condition.

We shall estimate the parameter 
$\theta$ by solving  the score equation 
$U_n(f, \theta) = o_P(n^{-1/2})$, where
$$
U_n(f, \theta) = {1 \over n} 
\sum_{i=1}^n \int_0^{\tau_0} \tilde W_f(u, \theta, Z_i)
 \tilde M_i(du, \theta)
$$
and 
$$
\tilde M_i(t, \theta)  =  N_i(t) - \int_0^t Y_i(u) \alpha
(\Gamma_{n\theta}(u), \theta, Z_i) \Gamma_{n\theta}(du) 
$$
Here 
$ \tilde W_f$ is defined by
substituting the estimate $\Gamma_{n\theta}$ into the 
score function 
$W_f(t,\theta, Z)$ of Proposition 2.1. 
A more explicit form of the score process is given
in Section 5. 

Let $\Sigma_0(f, \theta_0) = \Sigma_1(f, \theta_0) +
\Sigma_2(f, \theta_0)$ be given by (2.10)
and set 
$$
V(f, \theta_0) = \int_0^{\tau_0} \cov_0[f - \l' \phi_0, \dot \l + \l'  
D_0[\dot \l]](u) EN(du),
$$
where $\phi_0 = \phi^*$ is the solution to the Fredholm equation (2.6).
Note that if $f(x,\theta,Z) = \dot \l(x, \theta,Z)$, then
\beaa
V(f,\theta_0) & = & \Sigma_1(f, \theta_0) \\
& + & \int \cov_0[f - \l'\phi_0, \l'](u) (\phi_0 + D_0[\dot \l](u) EN(du)
\\
& = & \Sigma_1(f, \theta_0) + \Sigma_2(f, \theta_0)
\eeaa

\bcon
\bitem
\item[{(i)}] The matrix $\Sigma_0(f, \theta_0)$ is 
positive definite and the matrix 
 $V(f, \theta_0)$ is non-singular.
\item[{(ii)}]  
The estimate $\phi_{n\theta}$ of the solution to the Fredholm equation (2.6) satisfies
\newline $ \limsup_n \sup\{\Vert \phi_{n\theta} \Vertv: 
 \theta \in  B(\theta_0, \epsilon_n) \} = O_P(1)$ and
$ \sup\{\Vert \phi_{n\theta} - \phi_0 \Verti: 
 \theta \in  B(\theta_0, \epsilon_n) \} = o_P(1)$.
\eitem
\econ

The form of the solution $\phi_0$ to the equation (2.6) is given in
Section 5.3. Therein we also verify the condition 4.3 (ii) for the 
sample counterpart of this equation  based on an estimator
$\Gamma_{n\theta}$ satisfying the conditions 4.1.

\bprop Suppose that the conditions 2.1-2.3  and 4.1-4.3 hold.
 Then, with probability tending to 1, the score equation
$U_n(f, \theta ) = o_P(n^{-1/2}) $ has a solution 
in $B(\theta_0, \epsilon_n)$.
In addition, $\sqrt n(\hat \theta - \theta_0)$ converges
in distribution to a  multivariate normal 
variable $\N(0, \Sigma(f, \theta_0))$ 
with covariance matrix
$\Sigma(f, \theta_0) = 
(V^{-1} \Sigma_0[V^{-1}]^T)(f, \theta_0)$.
\eprop

In this proposition the asymptotic covariance matrix
of the estimate does not depend on the estimate of the 
unknown transformation. In addition, if we choose $f = \dot \l$, then
proposition 2.1 entails $\Sigma(f, \theta_0) = \Sigma_0(f, \theta_0)$.

The second version of this proposition, assumes that a
preliminary $\sqrt n$- consistent estimate $\hat \theta^{(0)}$ of
$\theta_0$ is available. Define
$$
\hat \theta = \hat \theta^{(0)} + V_n(f, \hat \theta^{(0)})^{-1}
U_n(f, \hat \theta^{(0)}) 
$$
Here $V_n(f, \hat \theta^{(0)})$ is the plug-in estimate
of the matrix $V_n(f, \theta_0)$. Section 5.2 gives the explicit form
of this matrix. 

\bprop
Suppose that the conditions 2.1-2.3 and 4.1-4.3 hold. Then 
$\sqrt n (\hat \theta - \theta_0)$ converges in distribution
to $\N(0, \Sigma(f, \theta_0))$.  
\eprop

Examples of simple $\sqrt n$ consistent estimators of the parameter $\theta$
were given in \cite{r5,r6,r7,r9,r12,r21}.

\sect{Proof of Proposition 3.1}

\subsect{An auxiliary  lemma}

The proof of Proposition 3.1 is based on the following modification
of Theorem 2 in Bickel \it et al. \rm [4, p.518]. 

\blem Let $B(\theta_0, \epsilon_n) = 
\{\theta: |\theta - \theta_0| \leq \epsilon_n)$ be a 
ball of radius $\epsilon_n, \epsilon_n \to 0$,
$\sqrt n \epsilon_n \to \infty$. Suppose that 
\bitem
\item[{(i)}] $\sqrt n U_n(\theta_0) \implies \N(0, \Sigma_0(\theta_0))$.
\item[{(ii)}] $V_{2n}(\theta_0) \to_P V(\theta_0)$. 
\item[{(iii)}] The matrices $\Sigma_0(\theta_0)$ and
$V(\theta_0)$ are nonsingular.
\item[{(iv)}] 
$U_n(\theta) - U_n(\theta_0) =
(\theta - \theta_0)^T V_n(\theta_0) + \remm(\theta)$,
where 
$$
\sup \left \{{\sqrt n |\remm(\theta) - \remm(\theta_0)| \over 1 + \sqrt n|\theta - \theta_0|}: \theta \in B(\theta_0,
\epsilon_n)  \right \} \to_P 0.
$$
\eitem
If the assumptions (i)-(iv) are satisfied then
with probability tending to 1, the score equation $U_n(\theta) =
o_P(n^{-1/2})$ has a solution $\hat \theta$ in $B(\theta_0, \epsilon_n)$
and 
$$
\sqrt n(\hat \theta -
\theta_0) \implies N(0, [V^{-1}\Sigma_0 (V^T)^{-1}](\theta_0))
$$
\elem

\Proof
Let
$\bar U_n(\theta) = U_n(\theta) - \remm(\theta)$. We have
 $\bar U_n(\theta_0) = U_n(\theta_0)$ because $\remm(\theta_0) = 0$. 
Set
\beaa
a_n & = & 
\Vert I - V^{-1}(\theta_0)V_n(\theta_0) \Vert = o_P(1), \\
A_n & = & V_2^{-1}(\theta_0) \bar U_n(\theta_0) = O_P(n^{-1/2}).
\eeaa
Finally, define 
$h_n(\theta) = \theta - V^{-1}(\theta_0) \bar U_n(\theta)$, and put
 $\theta^{(0)}_n = \theta_0$ and $\theta_n^{(m)} = h_n(\theta_n^{(m-1)})$
for $m \geq 1$.
The condition (iv) implies that for $m \geq 1$ we have
\bea
& & \theta_n^{(m)} - \theta_n^{(0)} =
[I - (V^{-1} V_{n})(\theta_0)] (\theta_n^{(m-1)} - \theta_n^{(0)})
- V^{-1}(\theta_0) \bar U_n(\theta_0) \\
& & h(\theta_n^{(m)}) - h(\theta_n^{(m-1)}) =
[I - (V^{-1} V_{n})(\theta_0)] (\theta_n^{(m)} - \theta_n^{(m-1)})
\eea

Similarly to Bickel et al ([4, p.518]), (5.1) - (5.2) implies that the
mapping $h_n$ is a contraction on the ball
$B_n = \{\theta: |\theta - \theta_0| \leq A_n/(1 - a_n)
\}$. With probability tending to 1, $B_n \subset B(\theta_0, \epsilon_n)$,
because $A_n = O_P(n^{-1/2}), a_n = o_P(1)$ and $\sqrt n \epsilon_n \uparrow
\infty$.  
It follows that with probability tending
to 1, the equation $V^{-1}(\theta_0)\bar U_n(\theta) = 0$ has a unique
root $\hat \theta$ in $B_n$ satisfying $\hat \theta - \theta_0 =
V_2(\theta_0)^{-1} \bar U_n(\theta_0) = O_P(n^{-1/2})$. We also have 
$U_n(\hat \theta) = \bar U(\hat \theta) + \remm(\hat \theta) =
o_P(|\hat \theta - \theta_0| +n^{-1/2}) = o_P(O_P(n^{-1/2}) +n^{-1/2}) = 
o_P(n^{-1/2})$.
 \done

\subsect{Proof of Proposition 3.1}

Define
\beaa
\hat s[1](t, \Gamma_{n\theta}, \theta) & = & {1 \over n}
\sum_{i=1}^n Y_i(t) \alpha(\Gamma_{n\theta}(t), \theta, Z_i), \\
\hat s[f](t, \Gamma_{n\theta}, \theta) & = & {1 \over n}
\sum_{i=1}^n Y_i(t) f (\Gamma_{n\theta}(t), \theta, Z_i) 
\alpha(\Gamma_{n\theta}(t), \theta, Z_i), \\
\hat e[f](t, \Gamma_{n\theta}, \theta) & = & 
{\hat s[f] \over \hat s[1]} (t, \Gamma_{n\theta}, \theta). 
\eeaa
Similarly to section 2, we put
\beaa
\hat \cov[f_1, f_2](t, \Gamma_{n\theta}, \theta) 
& = & (\hat e[f_1 f_2^T] - 
\hat e[f_1] \hat e[f^T_2])(t, \Gamma_{n\theta}, \theta), \\
\hat \var[f](t, \Gamma_{n\theta}, \theta) & = & 
\hat \cov[f,f](t, \Gamma_{n\theta}, \theta). 
\eeaa
Let $N_.(t) = n^{-1} \Sigma_{i=1}^n N_i(t)$ and set
\bea
\P_{n\theta}(s,t) & = & \exp - \int_s^t \hat s[\l']
(u, \Gamma_{n\theta}, \theta) C_{n\theta}(du), \nonumber \\
C_{n\theta}(t) & = & \int_0^t \hat s[1](u, \Gamma_{n\theta}, \theta)^{-2}
N_.(du) \\
\hat \rho[f, \phi_{n\theta}](t, \Gamma_{n\theta}, \theta)  
& = & \hat \cov[f, \l'](t, \Gamma_{n\theta}, \theta) -
\hat \var[\l'](t, \Gamma_{n\theta}, \theta)  \phi_{n\theta}(t),
\nonumber
\eea
where $ \phi_{n\theta}$ is an estimate of the solution to the 
Fredholm equation (2.6).

The score process for estimation of the parameter $\theta$ is given
by
$$
\tilde  U_n(f, \theta) = {1 \over n} \sum_{i=1}^n \int_0^{\tau_0}
\tilde  W_f(t, \theta, Z_i) \tilde  M_i(dt, \theta),
$$
where 
$$
\tilde  M_i(t, \theta) = 
N_i(t) - \int_0^t Y_i(u) \alpha(\Gamma_{n\theta}(u), \theta, Z_i) 
\Gamma_{n\theta}(du)
$$
and 
\beaa
\tilde  W_f(t, \theta, Z_i) & = &
b_{1i}(t, \Gamma_{n\theta}, \theta) - 
b_{2i}(t, \Gamma_{n\theta}, \theta)  \phi_{n\theta}(t), \\ 
& - & [\hat s[1](t, \Gamma_{n\theta}, \theta)]^{-1} 
\int_t^{\tau_0}  \P_{n\theta}(t,u) \hat \rho[f, \phi_{n\theta}](u,
\Gamma_{n\theta}, \theta) N_.(du), \\
b_{1i}(t, \Gamma_{n\theta},\theta) & = & 
f(\Gamma_{n\theta}(t), \theta, Z_i) -
\hat e[f](t, \Gamma_{n\theta}, \theta), \\   
b_{2i}(t, \Gamma_{n\theta},\theta) & = & 
\l'(\Gamma_{n\theta}(t), \theta, Z_i) -
\hat e[\l'](t, \Gamma_{n\theta}, \theta).    
\eeaa
The form of the score process simplifies if 
we introduce
$$
\hat \Gamma_{n\theta}(t) = 
\int_0^t {N_.(du) \over 
\hat s[1](u, \Gamma_{n\theta}, \theta)}.
$$
We have
\beaa
U_n(f, \theta) & = & {1 \over n} \sum_{i=1}^n
\int_0^{\tau_0} \left [ 
b_{1i}(t, \Gamma_{n\theta}, \theta) - 
b_{2i}(t, \Gamma_{n\theta}, \theta) \phi_{n\theta}(t)\right ]
N_i(dt) \\
& - & \int_0^{\tau_0} \left [
\int_t^{\tau_0} \P_{n\theta}(t,u) \hat \rho[f, \phi_{n\theta}](u,
\Gamma_{n\theta}, \theta) N_.(du) \right] [\hat \Gamma_{n\theta}
- \Gamma_{n\theta}](dt). 
\eeaa

Set $M_.(t) = n^{-1} \Sigma_{i=1}^n M_i(t)$. Then
\bea
[\hat \Gamma_{n0} - \Gamma_{n0}](t) & = & \int_0^t
{M_.(du) \over s[1](u, \Gamma_0(u), \theta_0)} -
[ \Gamma_{n0} - \Gamma_0](t) \nonumber \\
& - & \int_0^t [  \Gamma_{n0} - \Gamma_0](u) 
e_0[\l'](u)  \Gamma_0(du) + o_p(n^{-1/2}) \;.
\eea
The score process $U_n(f, \theta_0)$ can be represented as a 
sum $U_n(f, \theta_0) = \sum_{j=1}^4 U_{nj}(f, \theta_0)$, where
\beaa
U_{n1}(f, \theta_0) & = & {1 \over n} 
\sum_{i=1}^n \int_0^{\tau_0}
[ b_{1i}(t, \Gamma_0, \theta_0) - 
 b_{2i}(t, \Gamma_0, \theta_0)  \phi_{0}(t)] 
N_i(dt), \\
U_{n2}(f, \theta_0) & = &  - \int_0^{\tau_0}
[\hat \Gamma_{n0} - \Gamma_{n0}](du) \int_u^{\tau} 
\P_0(u,t)  \rho_0[f,\phi_0](t) EN_.(dt) + o_P(n^{-1/2}),\\
U_{n3}(f, \theta_0) & = &  - \int_0^{\tau_0} 
[\Gamma_{n0} - \Gamma_0](t) \rho_0[f,\phi_0](t) EN(dt) + o_P(n^{-1/2}),\\
U_{n4}(f, \theta_0) & = &  
 \int_0^{\tau_0} [\phi_{n0} - \phi_0](t)
{1 \over n} \sum_{i=1}^n b_{2i}(\Gamma_0(t), \theta_0,t) 
N_i(dt) = o_P(n^{-1/2})  \;.
\eeaa
By  central limit theorem, 
we have $\sqrt n U_{n1}(f, \theta_0) \implies
 N(0, \Sigma_1(f, \theta_0))$. The matrix $\Sigma_1(f, \theta_0)$ is defined in
 Section 2. Further, 
\beaa
& & U_{n2}(f, \theta_0) + U_{n3}(f, \theta_0) =
 -  \int_0^{\tau_0} 
\left [
\int_0^t [\hat \Gamma_{n0} - \Gamma_{n0}](du) \P_0(u,t)  \right]
\rho_0[f,\phi_0](t) EN(dt) \\
& - & \int_0^{\tau_0}
[ \Gamma_{n0} -  \Gamma_0](t) \rho_0[f,\phi_0](t)EN(dt) 
+ o_P(n^{-1/2}) = \\
& = & 
  \int_0^{\tau_0} 
\left [
\int_0^t [\hat \Gamma_{n0} - \Gamma_{n0}](du) \int_u^t e_0[\l'](s) 
\Gamma_0(ds) \P_0(s,t)  \right]
\rho_0[f,\phi_0](t) EN(dt) \\ 
& - & \int_0^{\tau_0} 
\left [
\int_0^t [\hat \Gamma_{n0} - \Gamma_{n0}](du)   \right]
\rho_0[f,\phi_0](t) EN(dt) \\
& - &  \int_0^{\tau_0}
[\Gamma_{n0} - \Gamma_0](t) \rho_0[f,\phi_0](t)EN(dt) 
+ o_P(n^{-1/2}) \\
& = & 
 \int_0^{\tau_0} 
\left [
\int_0^t [\hat \Gamma_{n0} - \Gamma_{n0}](s)  e_0[\l'](s) 
\Gamma_0(ds) \P_0(s,t)  \right]
\rho_0[f,\phi_0](t) EN(dt) \\
& - & \int_0^{\tau_0} 
[\hat \Gamma_{n0} - \Gamma_{n0}](t)
\rho_0[f,\phi_0](t) EN(dt) \\
& - &  \int_0^{\tau_0}
[ \Gamma_{n0} -  \Gamma_0](t) \rho_0[f,\phi_0](t)EN(dt) 
 + o_P(n^{-1/2}) \;. 
\eeaa
Next substitution of (5.3)  yields
\beaa
& & U_{n2}(f, \theta_0) + U_{n3}(f, \theta_0)  = \\
&- &  \int_0^{\tau_0} \left [ \int_0^t
\left (\int_0^s [\Gamma_{n0} - \Gamma_0](v) e_0[\l'](v) EN(dv) \right )
e_0[\l'](s) 
\Gamma_0(ds) \P_0(s,t)  \right]
\rho_0[f,\phi_0](t) EN(dt) \\
& - & \int_0^{\tau_0} 
\left [
\int_0^t  [  \Gamma_{n0} -  \Gamma_0](s)  e_0[\l'](s) 
\Gamma_0(ds) \P_0(s,t)  \right]
\rho_0[f,\phi_0](t) EN(dt) \\
& + & \int_0^{\tau_0} 
\left [
\int_0^t \left (\int_0^s {dM_. \over s_0[1]} \right)  e_0[\l'](s) 
\Gamma_0(ds) \P_0(s,t)  \right]
\rho_0[f,\phi_0](t) EN(dt) \\
& + & 
\int_0^{\tau_0} \left [\int_0^t [\Gamma_{n0}-\Gamma_0](u)e_0[\l'] EN(du)
\right] \rho_0[f,\phi_0](t)
EN(dt) \\
& - & \int_0^{\tau_0} \left (\int_0^t 
{dM_. \over s_0[1]} \right)  \rho_0[f,\phi_0](t) EN(dt)
+ o_P(n^{-1/2}) \;.
\eeaa
Using
\beq
\P_0(u,t) -1 = \int_u^t s_0[\l'](v) C_0(dv)\P_0(v,t) = \int_u^t e_0[\l'](v) \Gamma_0(dv) \P_0(v,t)
\eeq
and Fubini theorem, it is easy to see that the first, the second and 
the fourth term of this expansion sum to 0. The sum of the remaining
terms is
$$
U_{n2}(f, \theta_0) + U_{3n}(f, \theta_0) =
 - \int_0^{\tau_0} {M_.(du) \over s_0[1](u)} \left [\int_u^{\tau_0} \P_0(u,t)
\rho_0[f,\phi_0](t) EN(dt) \right] + o_P(n^{-1/2}) \;.
$$
We have 
$\sqrt n [U_{n2} + U_{n3}](f, \theta_0) \implies N(0, \Sigma_2(f, \theta_0))$,
and the matrix $\Sigma_2(f, \theta_0)$ is defined in Section 2.
It is also easy to verify that 
$\sqrt n [U_{n2} + U_{n3}](f,\theta_0)$ and $\sqrt n U_{n1}(f, \theta_0)$
are asymptotically uncorrelated. Therefore  $\sqrt n U_n(f, \theta_0)
\implies \N(0, \Sigma_0(f, \theta_0))$, $\Sigma_0 = \Sigma_1 + \Sigma_2$.

We consider now the expansion of the score process 
$U_n(f, \theta) - U_n(f, \theta_0)$ for 
$\theta \in B(\theta_0, \epsilon_n)$.
Set $\hat W(\theta) = \hat \Gamma_{n\theta} - \Gamma_{n\theta} - 
\hat \Gamma_{n0} + \Gamma_{n0}$. Then
\bea
\hat W(\theta)(t) & = & - (\theta - \theta') \int_0^t e_0[\dot \l](u)
\Gamma_0(du) -
\int_0^{t} [\Gamma_{n\theta} -
\Gamma_{n0}](u) e_0[\l'](u) \Gamma_0(du) \nonumber \\
& - & [\Gamma_{n\theta}
- \Gamma_{n0}](t) + \mbox{Rem}(\theta)(t),
\eea
where the remainder term satisfies
$$
\sup \{{\sqrt n|\mbox{Rem}(\theta)(t) - \mbox{Rem}(\theta_0)(t)| \over
\sqrt n|\theta - \theta_0| + 1}: \theta \not= \theta_0, \theta \in
B(\theta_0, \epsilon_n), t \leq \tau_0\} = o_P(1) 
$$
For $\theta \in B(\theta_0, \epsilon_n)$, we also have
$\Vert \Gamma_{n\theta} - \Gamma_{n0} \Verti = 
o_P(|\theta - \theta_0| + n^{-1/2})$.

Define
\beaa
\tilde e[f](u, \Gamma_0, \theta, \theta_0) &= &
{\Sigma_{i=1}^n Y_i(u) f(\Gamma_0(u), \theta, Z_i) \alpha(\Gamma_0(u), \theta_0, Z_i)
\over n \hat s[1](u, \Gamma_0, \theta_0)}, \\
\tilde e[\l'](u,\Gamma_0, \theta, \theta_0) & = &
{\Sigma_{i=1}^n Y_i(u) \l'(\Gamma_0(u), \theta, Z_i) \alpha(\Gamma_0(u), \theta_0, Z_i)
\over n \hat s[1](u, \Gamma_0, \theta_0)} 
\eeaa
ans let
\beaa
I_{1n}(\theta) & = & {1 \over n} \sum_{i=1}^n \int_0^{\tau_0}\biggl [
f(\Gamma_0(u),\theta, Z_i) - f( \Gamma_0(u), \theta_0, Z_i) \\
&  & \quad - \tilde e[f](u,\Gamma_0,\theta,\theta_0)
  + \hat e[f](u, \Gamma_0, \theta_0)\biggr ] N_i(du), \\
I_{2n}(\theta) & = &  -{1 \over n} \sum_{i=1}^n \int_0^{\tau_0} \biggl [
\l'(\Gamma_0(u), \theta,Z_i) - \l'(\Gamma_0(u),\theta_0 Z_i) \\
&  & \quad -\tilde e[\l'](u, \Gamma_0, \theta, \theta_0)  
 + \hat e[\l'](u,\Gamma_0,\theta_0)\biggr ] 
\phi_0(u) N_i(du).
\eeaa
The condition 4.2 implies that 
$I_{qn}(\theta), q = 1,2 $ is a mean zero square integrable martingale and 
$ \Var[\sqrt n I_{qn}(\theta)] = O(1) 
(\theta -\theta_0)^2 = O(\epsilon_n^2) = o(1)$.    
Using Hoeffding projection method, we can show that the right hand
side of these expressions can be approximated by a sum of U-processes
of degree $k, k \leq 4$ over  Euclidean classes of functions
with square integrable envelopes $H_{nk}, k \leq 4$. Under the
assumption of the transformation model, the $L_2$-norm of $EH_{nk}$ is
of order  $O(\epsilon_n)$. Application of maximal inequalities for
U-processes indexed by Euclidean classes of functions \cite{r24} shows
that $I_{qn}(\theta) = o_P(n^{-1/2})$, uniformly in $\theta \in B(\theta_0,\epsilon_n)$. The details are similar to \cite{r12,r25}, so we omit the proof.

The score process satisfies $U_n(f , \theta) - U_n(f, \theta_0)
= \Sigma_{j=1}^7 I_{jn}(\theta)$, where $I_{1n}(\theta)$ and 
$I_{2n}(\theta)$ are defined as above and
\beaa
I_{3n}(\theta) & = & {1 \over n}
\sum_{i=1}^n  \int_0^{\tau_0} [f(\Gamma_{n\theta}(u), \theta,Z_i) - f(\Gamma_{n0}(u), \theta_0, Z_i) \\
&  &  \quad - \hat e[f](u, \Gamma_{n\theta}, \theta) +
\hat e[f](u, \Gamma_{n0}, \theta_0)] N_i(du) - I_{1n}(\theta) \\
& = & {1 \over n}
\sum_{i=1}^n \int_0^{\tau_0} [ \Gamma_{n\theta} - \Gamma_{n0}](u)[f'(\Gamma_0(u), 
\theta_0, Z_i) - \hat e[f'](u,\Gamma_0, \theta_0)] N_i(du) \\      
& - & (\theta- \theta_0)\int_0^{\tau_0} \hat \cov [f, \dot \l]
(u, \Gamma_0, \theta_0)N_.(du) \\
& - & \int_0^{\tau_0}[\Gamma_{n\theta} - \Gamma_{n0}](u) 
\hat \cov[f, \l'](u, \Gamma_0, \theta_0) N_.(du) 
 +  
o_P(|\theta - \theta_0| + \Vert \Gamma_{n\theta} - \Gamma_{n0} \Verti + n^{-1/2}) \\
& = & -(\theta - \theta_0) \int_0^{\tau_0} \cov_0[f, \dot \l](u) N_.(du)
- \int_0^{\tau_0} [\Gamma_{n\theta} - \Gamma_{n0}] \cov_0[f,\l'](u) N_.(du) \\
& + &
o_p(|\theta - \theta_0| + n^{-1/2}) \\
I_{4n}(\theta) & = & 
-{1 \over n}
\sum_{i=1}^n  \int_0^{\tau_0} [\l'(\Gamma_{n\theta}(u), \theta,Z_i) - \l'(\Gamma_{n0}(u), \theta_0, Z_i) 
 \\ 
&  & \quad  -
\hat e[\l'](u, \Gamma_{n\theta}, \theta) +
\hat e[\l'](u, \Gamma_{n0}, \theta_0)]\phi_0(u) N_i(du) - I_{2n}(\theta) \\
& = & -{1 \over n}
\sum_{i=1}^n\int_0^{\tau_0} [ \Gamma_{n\theta} - \Gamma_{n0}](u)[\l''(\Gamma_0(u), 
\theta_0, Z_i) - \hat e[\l''](u,\Gamma_0, \theta_0)] \phi_0(u) N_i(du) \\      
& + & (\theta- \theta_0)\int_0^{\tau_0} \hat \cov [\l', \dot \l]
(u, \Gamma_0, \theta_0) \phi_0(u)N_.(du)  \\
& + & \int_0^{\tau_0}[\Gamma_{n\theta} - \Gamma_{n0}](u) 
\hat \var[l'](u, \Gamma_0, \theta_0) \phi_0(u)N_.(du) +
o_P(|\theta - \theta_0| + \Vert \Gamma_{n\theta} - \Gamma_{n0} \Verti + n^{-1/2}) \\
& = & (\theta - \theta_0) \int_0^{\tau_0} \cov_0[\l', \dot l](u) \phi_0(u)N_.(du)
+\int_0^{\tau_0} [\Gamma_{n\theta} - \Gamma_{n0}] \var_0[\l'](u) \phi_0(u)N_.(du) \\
& + &
o_p(|\theta - \theta_0| + n^{-1/2}) .
\eeaa
Combining,
\beaa
\sum_{j=1}^4 I_{jn}(\theta) & = & - (\theta - \theta_0)
\int_0^{\tau_0} \cov_0[f - \l' \phi_0, \dot \l](u) N_.(du) \\
& - & \int_0^{\tau_0} [\Gamma_{n\theta} - \Gamma_{n0}](u) \rho_0[f, \phi_0](u)
N_.(du) + o_P(|\theta - \theta_0| + n^{-1/2}).
\eeaa
The remaining three terms of the expansion are given by
\beaa
I_{n5}(f, \theta)  & = & 
- \int_0^{\tau_0} \hat W(\theta)(du) 
\int_u^{\tau_0}
\P_0(u,t)  \rho_0[f,\phi_0](u) N_.(du) \\
& = & (\theta - \theta_0)
\int_0^{\tau_0}  e_0[\dot \l] \Gamma_0(du) 
\int_u^{\tau_0}
\P_0(u,t)  \rho_0[f,\phi_0](u) N_.(du) 
\\
& + & \int_0^{\tau_0}  [ \Gamma_{n\theta} - \Gamma_{n0}](du) 
\int_u^{\tau_0} \P_0(u,t)  \rho_0[f,\phi_0](t) N(dt) \\
& + & \int_0^{\tau_0}  [\Gamma_{n\theta} - \Gamma_{n0}](u) 
e_0[\l'](u) \Gamma_0(du) 
\int_u^{\tau_0}  \P_0(u,t)  \rho_0[f,\phi_0](u) N(du) \\ 
& + & o_P(|\theta - \theta_0| + \Vert \Gamma_{n\theta} - 
\Gamma_{n0} \Vert
)\\
& = & - (\theta - \theta_0)\int_0^{\tau_0} D_0[\dot \l](u) \rho_0[f,\phi_0](u) 
N_.(du) \\
& + & \int_0^{\tau_0} [ \Gamma_{n\theta} - \Gamma_{n0}](u) 
\rho_0[f,\phi_0](f) N_.(du) 
 +  o_P(|\theta - \theta_0| + n^{-1/2}),
 \\
I_{n6}(f, \theta) & = & {1 \over n} \sum_{i=1}^n \int_0^{\tau_0}[ 
b_{2i}(u,\Gamma_{n\theta}(u), \theta) - b_{2i}(u, \Gamma_{n0},\theta_0)] 
[\phi_{n\theta} - \phi_0]N.(du) \\
& = & 
o_P(|\theta - \theta_0| + 
\Vert \Gamma_{n\theta} - \Gamma_{n0} \Vert
) = o_p(|\theta - \theta_0| + n^{-1/2}), \\
I_{n7}(f, \theta) & = & -\int_0^{\tau_0} \hat W(\theta)(du) \\
& & \times \quad \int_u^{\tau_0}[
 \P_{n\theta}(u,t) \hat \rho[f,\phi_{n\theta}](t, \Gamma_{n\theta}, \theta) - 
\P_0(u,t) \rho_0[f,\phi_0](t)] N_.(dt) \\
& - & \int_0^{\tau_0} [\hat \Gamma_{n0} - \Gamma_{n0}](du) \\
& & \times \quad \int_u^{\tau_0}
[\P_{n\theta}(u,t) \hat \rho[f,\phi_{n\theta}](t,\Gamma_{n\theta},\theta) - 
 \P_{n0}(u,t) \hat \rho[f,\phi_0](t, \Gamma_{n0},\theta_0)] N_.(dt) \\ 
& + & \int_0^{\tau_0} [\hat \Gamma_{n0} - \Gamma_{n0}](du) \\
& & \times \quad \int_u^{\tau_0}
[\P_{n0}(u,t) \hat \rho[f,\phi_{n0}](t, \Gamma_{n0},\theta_0) - 
 \P_0(u,t)  \rho_0[f,\phi_0](t)] N_.(dt) 
\\
& = & 
o_P(|\theta - \theta_0| + n^{-1/2}) 
 \;.
\eeaa
Using (5.4) and (5.5), we find that 
\beaa
U_n(f,\theta) - U_n(f, \theta_0) & = &
\Sigma_{j=1}^7
 I_{jn}(f, \theta) 
 = -(\theta - \theta_0) V_{n}(f, \theta_0) + 
 o_P(|\theta - \theta_0| + n^{-1/2}), \\
 V_{n}(f, \theta_0)  & = & 
\int_0^{\tau} \cov_0[f -\l' \phi_0, \dot \l +
 \l' D_0[\dot \l]](u) N_.(du). 
\eeaa
The matrix $V_{n}(f, \theta_0)$ converges in probability
to the matrix $V(\theta_0)$ defined in Section 4. 
so that Lemma 5.1 completes the proof.
\done

The proof of Proposition 4.2 follows from a similar expansion. The matrix
$\hat V_n(f, \theta^{(0)})$ can be defined by plugging-in the sample
counterpart of the covariance operator in the last display.

\subsect{Verification of the condition 4.3}

We have shown in \cite{r12} that the
 equation  (2.10) simplifies if we multiply both sides of it by 
$\P(0,t)^{-1}$. Let $
\tilde \psi(t)   =  \P(0,t)^{-1} \psi(t)$, 
$\tilde D[f] (t)   =   \P(0,t)^{-1} D[f](t)$ and
$\hat \rho[f, -D[f]](t)  =   \P(0,t)
 \rho[f,-D[f]](t)$
Set
$$
c(t)   =  \int_0^t \P(0,u)^{-2} dC(u), \quad
b(t)     = \int_0^t
\P(0,u)^2 B(du). 
$$
Multiplication of  (2.7) by
$\P(0,t)^{-1}$ yields  
\beq
\tilde \psi(t) + \int_0^\tau k(t,u) \tilde \psi(u) b(du) =
\int_0^{\tau} k(t,u) \tilde  \rho[f, - D[f]](u) EN(du) \;,
\eeq
where the  kernel $k$ is given by 
$
k(t,u) = c(t \wedge u).
$
The square integrability condition 2.3 is equivalent to
the assumption that $\kappa(\tau_0) = \int_0^{\tau_0} c(u) b(du)$ is
finite.
The solution to the equation (5.6) is given by
\beq
\tilde \psi(t) = \int_0^{\tau} \tilde \Delta(t,u) 
\tilde \rho[f, -D[f](u) EN(du) \;,
\eeq
where  $\tilde \Delta(t,u) = \tilde \Delta(t,u, -1)$, 
and $\tilde \Delta(t,u,\lambda)$ is the resolvent corresponding to the 
kernel $k$. 
 The solution to the equation (2.6) is given by
\beq
\phi(t) = - D[f](t) + \int_0^{\tau_0} \tilde \Delta(t,u) 
\rho[f,-D[f]](u) EN(du)\P(0,u) \P(0,t) 
\eeq

From \cite{r12},  $\tilde \Delta(u,t)$ is given by
\beq
\tilde \Delta(u,t) = {\Psi_1(0, u \wedge t) \Psi_0(u \vee t, \tau_0) \over
\Psi_0(0,  \tau_0)},
\eeq
where for $s < t$, the interval functions $\Psi_0(s,t)$ and $\Psi_1(s,t)$
are defined as solutions to the Volterra equations
\beaa
\Psi_0(s,t) &= & 1 + \int_{(s,t]} c((s,u])b(du) \Psi_0(u,t) =
1 + \int_{(s,t]} \Psi_0(s,u-)c(du) b([u,t]) \\
\Psi_1(s,t) & = & c((s,t]) + \int_{(s,t]} c((s,u])b(du) \Psi_1(u,t) = 
c((s,t]) + \int_{(s,t]} \Psi_1(s,u)b(du)c((u,t]).
\eeaa
Define also
\beaa
\Psi_2(s,t) &= & 1 + \int_{[s,t)}b([s,u)) c(du) \Psi_2(u,t) = 1 + 
 \int_{[s,t)} \Psi_2(s,u+)b(du) c((u,t)) \\
\Psi_3(s,t) & = & b([s,t)) + \int_{[s,t)} b([s,u))c(du) \Psi_3(u,t) = 
b([s,t)) + \int_{[s,t)} \Psi_3(s,u)c(du)b([u,t)) \;
\eeaa
Then
\beaa
\Psi_0(s,t) & = & 1 + \int_{(s,t]} \Psi_1(s,u) b(du) = 1 + \int_{(s,t]} c(du)
\Psi_3(u,t+) \:,\\
\Psi_1(s,t) & = & \int_{(s,t]} \Psi_0(s,u-) c(du) =  \int_{(s,t]} c(du)
\Psi_2(u,t+) \;,\\
\Psi_2(s,t) & = & 1 + \int_{[s,t)}  b(du) \Psi_1(u,t-) = 
1 + \int_{[s,t)} \Psi_3(s,u) c(du) \;, \\
\Psi_3(s,t) & = & \int_{[s,t)} \Psi_2(s,u)  b(du) =  
\int_{[s,t)}  b(du) \Psi_0(u,t-) \;.
\eeaa

If $\tau_0$ is an atom of the survival function $P(X > t)$, then
$\Psi_j, j = 0,1,2,3$ form bounded monotone increasing interval functions.
In particular, $\Psi_0(s,t) \leq \exp \kappa(\tau_0)$ and $\Psi_1(s,t)
\leq \Psi_0(s,t) [c(t)-c(s)]$. 
If $\tau_0$ is a continuity point of the survival function
$P(X > t)$ and $\kappa(\tau_0) < \infty$, then $\Psi_0(s, t) 
\leq \exp \kappa(\tau_0)$ for any $0 < s < t \leq \tau_0$, while
the remaining functions are locally bounded \cite{r12}.  

Next suppose that the transformation model holds with
$(\theta, \Gamma) = (\theta_0, \Gamma_0)$.
We assume  that the estimate $\Gamma_{n\theta}$ satisfies the 
conditions 4.1 and show the natural plug-in estimator $\phi_{n\theta}$
of (5.9) satisfies the conditions 4.3.
Define 
\begin{eqnarray*}
c_{n\theta}(t) & = & \int_0^t  {\cal P}_{n\theta}(0,u)^{-2} C_{n\theta}(du)\\
b_{n\theta}(t) & = &
\int_0^t {\cal P}_{n\theta}(0,u)^2 
\hat \var[\l'](u, \Gamma_{n\theta}, \theta) N_.(du),
\end{eqnarray*}
where $C_{n\theta}$ and $\P_{n\theta}$ are defined as in Section 5.2.
Then the sample analogue of the equation (5.9)
reduces to a system of linear equations which can be solved by
inverting a bandsymmetric tridiagonal matrix \cite {r11,r12}.

Denote by $\Psi_{n\theta,j}$ the sample counterparts of the 
interval functions $\Psi_j$, $j = 1,2,3,4$.
Using Fubini theorem we can show that 
\beaa
[ \Psi_{n\theta,0} - \Psi_0](s,t) & = & 
\int_{(s,t]} \Psi_{n\theta,0}(s,u-)c_{n\theta}(du) 
[b_{n\theta} - b]([u,t])\\  
& + & \int_{(s,t]} [c_{n\theta} - c)]((s,u]) b(du) \Psi_0(u,t) \\ 
& + & \int_{s < u_1 < u_2 \leq t}  \Psi_{n\theta,0}(s,u_1-)c_{n\theta}(du_1) 
[b_{n\theta} - b]([u_1,u_2)) c(du_2) \Psi_3(u_2, t+) \\ 
& + & \int_{s < u_1 < u_2 \leq t}  \Psi_{n\theta 1}(s,u_1)b_{n\theta}(du_1) 
[c_{n\theta} - c]((u_1,u_2]) b(du_2) \Psi_0(u_2, t) 
\eeaa
and
\beaa
[\Psi_{n\theta,1} - \Psi_1](s,t) & = & 
[c_{n\theta} - c]((s,t]) \\
& + & \int_{(s,t]} \Psi_{n\theta,1}(s,u) b_{n\theta}(du) 
[c_{n\theta} - c]((u,t])\\  
& + & \int_{(s,t]} [ \Psi_{n\theta,0} - \Psi_0](s,u-) c(du)
\eeaa
 
Under assumptions of the condition 4.1, we have
$
c_n  =  \sup \{|c_{n\theta} - c|(t): \theta \in B(\theta_0, \epsilon_n),
t \leq \tau_0 \} \to_P 0$ and
$b_n  =  \sup\{|b_{n\theta} - b|(t\pm): \theta \in B(\theta_0, \epsilon_n),
t \leq \tau_0\} \to_P 0$. In addition, for $q = 0,1$, we have
\beaa
& & \limsup_n \sup \{\Psi_{n\theta,q}(0,\tau): \theta \in B(\theta_0, \epsilon_n) \} \\
& & \quad \quad \leq [c(\tau_0) - c(0)]^q \exp[\kappa(\tau_0)](1 + o_p(1)) = O_p(1)
\eeaa  
Hence 
\beaa
\hat r &= & 
\sup\{|\Psi_{n\theta,0} - \Psi_0|(s,t): 0 < s < t \leq \tau, \theta
\in B(\theta_0, \epsilon_n) \}  \\
 &\leq &
2 b_n  \Psi_{n\theta,0}(0, \tau_0) \Psi_0(0,\tau_0) 
 +  2c_n \Psi_{n\theta,1}(0, \tau) \Psi_3(0,\tau) \\
& & \sup\{|\Psi_{n\theta,1} - \Psi_1|(s,t): 0 < s < t \leq \tau, \theta
\in B(\theta_0, \epsilon_n) \}  \\ 
&\leq &
2 c_n \Psi_{n\theta,0}(0, \tau_0) + 
\hat r c(\tau_0)
\eeaa
and both terms converge in probability to 0. 

Let $\tilde \Delta_{n\theta}(s,t)$ be defined similarly to (5.10)
The preceding calculations, imply that 
$\sup\{\tilde \Delta_{n\theta}(s,t): s,t \in [0,\tau_0], \theta \in
B(\theta_0,\epsilon_n) \} = O_P(1)$ and
$\sup\{|\tilde \Delta_{n\theta} - \tilde \Delta|(s,t): s,t \in [0, \tau_0],
\theta \in B(\theta_0, \epsilon_n)\} = o_P(1)$.
Verification that the sample analogue $\phi_{n\theta}$ of the equation
(5.9) satisfies the conditions 4.3 can be completed  using Gronwall's 
inequalities given in \cite{r12} and integration by parts.

\end{document}